\documentclass{article}
\usepackage{amsmath}
\usepackage{amsthm}
\newtheorem{thm}{Theorem}
\newtheorem{lem}{Lemma}
\usepackage[boxed]{algorithm2e}

\newtheorem{qn}{Question}

\makeatletter
\newcommand\ackname{Acknowledgements}
\if@titlepage
  \newenvironment{acknowledgements}{%
      \titlepage
      \null\vfil
      \@beginparpenalty\@lowpenalty
      \begin{center}%
        \bfseries \ackname
        \@endparpenalty\@M
      \end{center}}%
     {\par\vfil\null\endtitlepage}
\else
  \newenvironment{acknowledgements}{%
      \if@twocolumn
        \section*{\abstractname}%
      \else
        \small
        \begin{center}%
          {\bfseries \ackname\vspace{-.5em}\vspace{\z@}}%
        \end{center}%
        \quotation
      \fi}
      {\if@twocolumn\else\endquotation\fi}
\fi
\makeatother

\begin{document}

\title{On the sum of two squares and at most two powers of $2$}
\markright{Sum of squares and two powers of $2$}
\author{David J. Platt and Timothy S. Trudgian}

\maketitle

\begin{abstract}
We demonstrate that there are infinitely many integers that cannot be expressed as the sum of two squares of integers and up to two non-negative integer powers of $2$.
\end{abstract}

\section{Introduction}

%In what follows, squares will always refer to non-zero squares of integers, and to non-negative integer powers of two.

It is well known that that are infinitely many integers that cannot be expressed as the sum of two squares, to whit any whose prime factorisation contains a prime $p\equiv 3 \pmod{4}$ to an odd power --- see \cite[Thm 278]{HardyW}.
Almost as trivially, there are infinitely many integers that cannot be expressed as the sum of two squares and at most one power of $2$ --- see Theorem \ref{thm:one2} below.

Crocker \cite{Crocker2008} proved that one can generate an infinitude of integers not expressible as a sum of two squares and at most two powers of $2$ provided one can show the existence of a single integer $N_0 \equiv 0\pmod{36}$ that cannot be so expressed\footnote{In fact the condition $0\pmod{18}$ will suffice as we show in Lemma \ref{lem:tower}.}. Crocker lists $142$ congruence conditions on such an $N_0$ and proves that the first example must be below $2^{1417} = 3.62\ldots \times 10^{426}$. We give a much shorter proof of this and show
\begin{thm}\label{spade}
The smallest integer, greater than $1$, which cannot be represented as a sum of two squares and at most two powers of $2$ is $535\,903$. Moreover, for any $\alpha\geq 0$, no integer of the form $2^{\alpha} \times 1\,151\,121\,374\,334$ can be so expressed.
\end{thm} 

We first tackle the sums of squares and one power of $2$ in Section \ref{trowel}, which allows us to deal with two powers of $2$ in Section \ref{grout}. We give some details on our computations that allow us to prove Theorem \ref{spade} in Section \ref{sod}. Finally, we pose some open problems in Section \ref{clay}.

Throughout this paper `squares' denotes squares of integers, and `powers' denotes non-negative integral powers.

%We will first dispense with the one power of two case before considering the two powers of two situation.

\section{One power of $2$}\label{trowel}

We start with a simple Lemma that will be used in both the one and two powers of $2$ cases.

\begin{lem}\label{lem:times2}
Suppose an integer $n$ cannot be expressed as the sum of two squares. Then neither can $2^a n$ for $a\geq 0$.
\begin{proof}
Since $n$ cannot be expressed as the sum of two squares, its prime factorisation must contain a prime $p\equiv 3\pmod{4}$ to an odd power. This remains true after multiplying by any power of $2$.
\end{proof}
\end{lem}

We now focus one the one power of $2$ case.

\begin{lem}\label{lem:tower_one2}
Suppose we have an even positive integer $n$ which cannot be expressed as the sum of two squares and at most one power of $2$. Then neither can $2^\alpha n$ for any $\alpha\geq 0$.
\begin{proof}
Since neither $n$ nor $n-2^a$ with $a\geq 0$ can be expressed as the sum of two squares, then by Lemma \ref{lem:times2} we can say the same for $2n$ and $2n-2^{a+1}$. This leaves $2n-1$ which is $\equiv 3 \pmod{4}$ establishing the result for $2n$. The Lemma now follows by induction.
\end{proof}
\end{lem}

It is now trivial to find the first such $n$.

\begin{thm}\label{thm:one2}
The are infinitely many integers that cannot be expressed as the sum of two squares and at most one power of $2$.
\begin{proof}
Referring to \cite{OEIS_A274050}, we see that $142$ cannot be expressed\footnote{$142=2\cdot \mathbf{71}$, $141=\mathbf{3}\cdot\mathbf{47}$, $140=2^2\cdot5\cdot\mathbf{7}$, $138=2\cdot\mathbf{3}\cdot\mathbf{23}$, $134=2\cdot\mathbf{67}$, $126=2\cdot 3^2\cdot\mathbf{7}$, $110=2\cdot 5\cdot\mathbf{11}$, $78=2\cdot\mathbf{3}\cdot 13$ and $14=2\cdot\mathbf{7}$.} as the sum of two squares and at most one power of $2$. Theorem \ref{thm:one2} now follows\footnote{In fact we could start our ``tower'' at $71$.}
 via Lemma \ref{lem:tower_one2}.\end{proof}
\end{thm}

\section{Two powers of $2$}\label{grout}

For two powers of $2$, we proceed along similar lines.

\begin{lem}\label{lem:tower}
Suppose an integer $n\equiv 0 \pmod{18}$ cannot be expressed as the sum of two squares and at most two powers of $2$. Then neither can $2^\alpha n$ for any integer $\alpha\geq 0$.
\begin{proof}
Since none of $n$, $n-2^a$ with $a\geq 0$ nor $n-2^a-2^b$ with $a,b\geq 0$ are the sum of two squares, then by Lemma \ref{lem:times2} we can dispense with $2n$, $2n-2^{a+1}$ and $2n-2^{a+1}-2^{b+1}$.

This leaves $2n-1$ and $2n-1-2^b$ with $b\geq 1$, since the case $2n - 1 - 1$ is covered by Lemma \ref{lem:times2}. Now since $2n\equiv 0\pmod{4}$ we have immediately that neither $2n-1$ nor $2n-1-2^b$ with $b\geq 2$ can be expressed as the sum of two squares. 

Consider now $2n-3\equiv 6\pmod9$. The only values of $x^{2} \pmod9$ are $0$, $1$, $4$ and $7$. Therefore $2n-3$ cannot be the sum of two squares, whence the Lemma follows by induction.
\end{proof}
\end{lem}

\begin{lem}\label{lem:N0}
The number $N_0=1\,151\,121\,374\,334$ cannot be written as the sum of two squares and at most two powers of $2$.
\begin{proof}
We factorise $N_0$, $x=N_0-2^a$ and $y=N_0-2^a-2^b$ with $x,y,a,b\geq 0$ and $a>b$. In each case (and since $2^{40}+2^{35}<N_0<2^{40}+2^{36}$ there are $858$ of them to check), there is at least one prime $p\equiv 3\mod\; 4$ appearing to an odd power.
\end{proof}
\end{lem}

\begin{thm}
There are infinitely many integers that cannot be expressed as the sum of two squares  and at most two powers of $2$.
\begin{proof}
Since $18|N_0$ this is a simple corollary of Lemmas \ref{lem:tower} and \ref{lem:N0}.
\end{proof}
\end{thm}

\section{Computational aspects}\label{sod}

%In the one power of two case, the first example of an integer in the right form was $142$ which entailed little or no computational effort on our part. However, our 
Finding $N_0$ for the two powers of $2$ case required some ingenuity. We proceeded by implementing a simple sieve in ``C++'' (see Algorithm \ref{alg:sieve}). 

\begin{algorithm}\label{alg:sieve}
\SetAlgoLined
$S\leftarrow $ length of sieve\;
\For{$n\leftarrow 1$ \KwTo $\infty$}{
   \lIf{$n^2\geq S$}{{\bf break}}
   \For{$m\leftarrow 1$ \KwTo $n$}{
      \lIf{$n^2+m^2>S$}{{\bf break}}
      \lIf{$18|n^2+m^2$}{cross out $n^2+m^2$}
      \For{$a\leftarrow 0$ \KwTo $\infty$}{
         \lIf{$n^2+m^2+2^a>S$}{{\bf break}}
         \lIf{$18|n^2+m^2+2^a$}{cross out $n^2+m^2+2^a$}
         \For{$b\leftarrow 0$ \KwTo $a-1$}{
            \lIf{$n^2+m^2+2^a+2^b>S$}{{\bf break}}
            \lIf{$18|n^2+m^2+2^a+2^b$}{cross out $n^2+m^2+2^a+2^b$}}}}}
\caption{A Simple Sieve}
\end{algorithm}

Each integer 
%in the sieve range 
divisible by $18$ was represented by a single byte in a vector and initially all bytes were set to $1$. When a way of representing such an integer as a sum of two squares and at most two powers of $2$ was found, then the relevant entry in the vector was set to $0$. Crucially, because no reads of this vector were required and every possible write was of a $0$, we could allow several concurrent processes to perform the sieve on a block of shared memory in parallel without the need for memory locks. This was achieved using the Posix Threads library. Note that had we used each byte to represent $8$ integers and achieved crossing out via bit twiddling, we would no longer have been able to ignore potential memory clashes and performance would have suffered accordingly.

We used a single node of the University of Bristol's BlueCrystal Phase III cluster \cite{ACRC2015} which consists of $16\times 2.6$GHz Sandy Bridge cores sharing $256$Gbytes of memory. We set the length of the sieve to be $18\times 2^{36}$ so that the sieve vector occupied $64$ Gbytes and used $32$ threads. The elapsed time was a little over $31\frac{1}{2}$ hours and only the one suitable $N_0$ was found.

\section{Some open problems}\label{clay}
There is a rich history in representing numbers as the sum of one prime and powers of $2$. We briefly outline two problems below, and present questions on their parallels with sums of two squares and powers of $2$.

Romanov \cite{Romanov} proved that there is a positive proportion of integers that can be written as the sum of a prime and one power of $2$; van der Corput \cite{Corpus} proved there is a positive proportion of integers that cannot be so written. 

It is easy to see that no integer of the form $23 \pmod{72}$ can be written as a sum of squares and at most one power of $2$. Note first that $23 + 72k \equiv 3 \pmod 4$ and $23 + 72k - 2^{\alpha} \equiv 3 \pmod 4$ for all $\alpha \geq 2$, whence neither is a sum of two squares. All one needs to show now is that neither $23 + 72k - 1$ nor $23 + 72k -2$ is a sum of squares. The first is $6 \pmod8$ and the second is $3\pmod 9$, neither of which can be written as a sum of two squares. Thus $1/72 = 1.38\ldots \%$ of integers cannot be written as a sum of two squares and one power of $2$. 

\begin{qn}
Can one obtain good estimates on the density of numbers that can, and cannot, be written as a sum of two squares and one power of $2$?
\end{qn}

Crocker \cite{Pacific} proved that there are infinitely many odd $n$ not of the form $p + 2^a + 2^b$. Pan \cite{Pan} proved that the number of such $n$ that are less than $N$ is $\gg N^{1-\epsilon}$.
We compare these results to the set of numbers generated in Theorem \ref{spade}: there are $\gg \log N$ such integers less than $N$. Is this density close to the mark?

\begin{qn}
Can one obtain a good estimate on the density of integers that are not sums of two squares and at most two powers of $2$?
\end{qn}

Finally, we ask the following

\begin{qn}
Does there exist a constant $C$ such that every sufficiently large integer $N$ can be expressed as a sum of two squares and at most $C$ powers of $2$?
\end{qn}
If such a $C$ exists then clearly $C\geq 3$. Does $C=3$?  For two powers of $2$ we needed only to consider congruences modulo $2$ and modulo $9$ to prove Lemma \ref{lem:tower}. In principle one could try to develop a system of congruences to tackle the $C=3$ case --- though this would be a formidable operation!

We conclude by observing that if we can write $N>1$ as the sum of two squares and at most two powers of $2$, then we can write $N+2^\alpha$ as the sum of two squares and at most three powers of $2$. There are $123\,494$ $N\in[2,2^{36}]$ that cannot be expressed as the sum of two squares and two powers of $2$, but in every case, $N-2$ can be. Thus we can state that if there is an $N>1$ that cannot be written as the sum of two squares and at most three powers of $2$, then $N>2^{36} = 6.8\ldots \times 10^{10}$. This calculation took 34 hours on 32 cores of 2.8GHz AMD Opteron(tm) 6320.

\begin{acknowledgements}
The authors are grateful to Roger Heath-Brown for pointing out this problem and for his helpful suggestions. 
\end{acknowledgements}

%\begin{biog}
%\item[Woodrow Wilson] (twwilson@princeton.edu) received his PhD in history and political science from Johns Hopkins University. He held visiting positions at Cornell and Wesleyan before joining the faculty at Princeton, where he was eventually appointed president of the university.  Among his proudest accomplishments was the abolition of eating clubs at Princeton on the grounds that they were elitist.
%\begin{affil}
%Office of the President, Princeton University, Princeton NJ 08544\\twoodwilson@princeton.edu
%\end{affil}

%\item[Herbert Hoover] (hchoover@stanford.edu) entered Stanford University in 1891, after failing all of the entrance exams except mathematics.  He received his BS degree in geology in 1895, spent time as a mining engineer, then was appointed by his co-author to the U.S. Food Administration and the Supreme Economic Council, where he orchestrated the greatest famine relief efforts of all time.
%\begin{affil}Hoover Institution, Stanford University, Stanford CA 94305\\ herbhoover@stanford.edu
%\end{affil}
%\end{biog}
%\vfill\eject

\end{document}